\documentclass[hidelinks,final]{siamart190516}
\newcommand{\linelabel}[1]{}
\usepackage{amsmath,amssymb,mathtools}
\usepackage{smalljeff}
\usepackage{mathdots}
\newcommand{\E}{\mathbb{E}}
\newcommand{\prob}{\mathbb{P}}

\usepgfplotslibrary{groupplots}
\usepgfplotslibrary{fillbetween}
\DeclareMathOperator{\erf}{erf}

\DeclareMathOperator{\rank}{rank}

\newcommand{\tikzmark}[1]{\tikz[overlay,remember picture] \node (#1) {};}
\newcommand{\DrawBox}[1][]{    \tikz[overlay,remember picture]{
    \draw[black,#1]
      ($(left)+(-0.2em,0.9em)$) rectangle
      ($(right)+(0.2em,-0.3em)$);}
}

\newcommand{\TheTitle}{A Data-Driven McMillan Degree Lower Bound}
\newcommand{\TheAuthors}{Jeffrey M. Hokanson}
\headers{\TheTitle}{\TheAuthors}

\title{{\TheTitle}\thanks{Submitted to the editors 17 August 2018.
\funding{This work was supported by NSF VIGRE grants DMS-0240058 and DMS-0739420
	and the DARPA program Enabling Quantification of Uncertainty in Physical Systems (EQUiPS).
}}}

\author{
  Jeffrey M. Hokanson\thanks{
	Department of Computer Science, University of Colorado Boulder,
	1111 Engineering Dr, Boulder, CO 80309,
	(\email{Jeffrey.Hokanson@colorado.edu}, \url{http://www.hokanson.us}).}
}

\begin{document}
\maketitle
\begin{abstract}
In the context of linear time-invariant systems,
the McMillan degree prescribes the smallest possible dimension
of a system that reproduces the observed dynamics.
When these observations take the form of impulse response measurements
where the system evolves without input from an unknown initial condition,
a result of Ho and Kalman reveals the McMillan degree 
as the rank of a Hankel matrix built from these measurements.
Unfortunately using this result in experimental practice is challenging
as measurements are invariably contaminated by noise
and hence the Hankel matrix will almost surely be full rank.
Hence practitioners estimate the rank of this matrix---and thus the McMillan degree---by manually setting a threshold separating  
large singular values corresponding to the non-zero singular values of the noise-free Hankel matrix
and small singular values corresponding to perturbation of zero singular values of the noise-free Hankel matrix. 
Here we replace this manual threshold 
with a threshold guided by Weyl's theorem.
Specifically, assuming measurements are perturbed by additive Gaussian noise
we construct a probabilistic upper bound
on how much the singular values of the noise-free Hankel matrix can be perturbed;
this provides a conservative threshold for estimating the rank, and hence the McMillan degree.
This result follows from a new probabilistic bound on the 2-norm 
of a random Hankel matrix with normally distributed entries. 
Unlike existing results for random Hankel matrices, 
this bound features no unknown constants
and, moreover, is within a small factor of the empirically observed bound.
This bound on the McMillan degree provides an inexpensive alternative to more general
model order selection techniques such as the Akaike Information Criteria (AIC).
\end{abstract}
 
\section{Introduction}
Here we consider discrete-time, linear time-invariant dynamical systems that map an input $\ve u \in \ell_\infty(\mathbb{N})$ to an output 
$\ve y \in \ell_\infty(\mathbb{N})$.
Such systems are uniquely defined via their impulse response $\ve h \in \ell_1(\mathbb{N})$
though a discrete convolution
\begin{equation}
	\ve y = \ve h * \ve u \quad \text{where} \quad 
	[\ve h * \ve u]_k \coloneqq \sum_{j=0}^k h_j u_{k-j}.
\end{equation} 
In \emph{system identification}~\cite{GP77},
the goal is to recover the system described by $\ve h$
through observations of pairs of inputs $\ve u$ and outputs $\ve y$.
Rather than recovering $\ve h$ explicitly,
typically one recovers a \emph{state-space model} instead.
State-space models take the form 
\begin{equation}\label{eq:ss}
	\left\lbrace \begin{split}
		\ve x_{j} &= \ma A \ve x_{j-1} + \ve b u_j, \quad \ve x_{\,-1} = \ve 0 \\
		y_{j} &= \ve c^*\ve x_{j} 
	\end{split}	
	\right\rbrace
	\quad 
	\text{where} \quad \ve x_j, \ve b, \ve c \in \C^q,
	\ \ma A \in \C^{q\times q}.
\end{equation}
An important hyperparameter for many system identification algorithms 
is the dimension the state-space $q$ in~\cref{eq:ss}.
However, this dimension is not unique.
For example, the impulse response of 
the system in~\cref{eq:ss} is $\lbrace \ve c^* \ma A^k \ve b \rbrace_{k=0}^\infty$;
an enlarged system with 
$\ma A \to \begin{bsmallmatrix} \ma A & \ma 0 \\ \ma 0 & \star \end{bsmallmatrix}$,
$\ve b \to \begin{bsmallmatrix} \ve b \\ \ma 0 \end{bsmallmatrix}$,
and $\ve c \to \begin{bsmallmatrix} \ve c \\ \ve 0 \end{bsmallmatrix}$
has the same impulse response.
As identifying smaller systems requires fewer parameters and less computation,
we ask: what is the smallest possible state-space system 
whose impulse response is~$\ve h$?
This is the \emph{McMillan degree}~\cite[Rmk.~6.7.4]{Son98},
\begin{equation}\label{eq:mcmillan_def}
	\op M(\ve h) \coloneqq
		\min_{q\in \mathbb{N}} q \ 
			\text{ s.t. } \exists \  \ve b, \ve c \in \C^q, \ \ma A \in \C^{q\times q} 
			\text{ with } h_k = \ve c^* \ma A^k \ve b
			\ \  \forall k \in \mathbb{N},
\end{equation}
named in honor of Brockway McMillan's pioneering work on this subject~\cite{Mil52a,Mil52b}. 
Remarkably, the McMillan degree can be computed without
explicitly recovering a \emph{minimal realization} with matrices $\ma A$, $\ve b$, and $\ve c$
appearing in the optimizer of~\cref{eq:mcmillan_def}.
\begin{theorem}[Ho and Kalman {\cite[Thm.~2 Cor.]{HK66}}]\label{thm:HK}
Let $\ma H_\infty$ denote the infinite Hankel matrix built from~$\ve h\in \ell_1(\mathbb{N})$,
then 
\begin{equation}
	\op M(\ve h) = \rank(\ma H_\infty) \coloneqq \sup_{n \in \mathbb{N}} \ \rank(\ma H_{n,n})
	\quad \text{where} \quad
	\ma H_\infty \coloneqq \begin{bmatrix}
		h_0 & h_1 & h_2 &\cdots \\
		h_1 & h_2 & h_3 &\cdots \\
		h_2 & h_3 & h_4 & \cdots \\
		\vdots & \vdots & \vdots & \ddots 
	\end{bmatrix}
\end{equation}
and $\ma H_{m,n} \in \C^{m\times n}$ denotes the $m\times n$ leading principal submatrix of $\ma H_\infty$. 
\end{theorem}

At first glance, it would appear straightforward to apply 
Ho and Kalman's result to estimate the McMillan degree in experimental practice.
By providing an impulse response input, $\ve u = \lbrace 1,0,0,0,\ldots\rbrace$,
we can directly measure the impulse response $\ve h$ in the output $\ve y$.
However, this poses two challenges.
One is that we necessarily only measure finite data
and hence cannot build the infinite-dimensional matrix $\ma H_\infty$.
Instead we can only construct a lower bound on the McMillan degree
from the rank of $\ma H_{m,n}$.
A more substantive challenge is that experimental measurements 
are invariably contaminated with noise.
So rather than measuring $h_k$,
we can only measure a noisy version $\widetilde{h}_k$.
If we build the analogous Hankel matrix $\tma H_{m,n}\in \C^{m\times n}$ from $\widetilde{h}_k$,
\begin{equation}
	\setlength\arraycolsep{3pt}
	\ma H_{m,n}
	\coloneqq \begin{bmatrix}
			 {h}_0 &  {h}_1 & \cdots &   {h}_{n-1} \\
			 {h}_1 &  {h}_2 & \cdots & {h}_{n} \\
				\vdots & & & \vdots \\
			 {h}_{m - 1} &  {h}_{m} & \cdots &  {h}_{m+n-2}
		\end{bmatrix}
	\quad
	\tma H_{m,n}
	\coloneqq \begin{bmatrix}
			 \widetilde{h}_0 &  \widetilde{h}_1 & \cdots &   \widetilde{h}_{n-1} \\
			 \widetilde{h}_1 &  \widetilde{h}_2 & \cdots & \widetilde{h}_{n} \\
				\vdots & & & \vdots \\
			 \widetilde{h}_{m - 1} &  \widetilde{h}_{m} & \cdots &  \widetilde{h}_{m+n-2}
		\end{bmatrix},
\end{equation}
then even if $\ma H_{m,n}$ is low rank, $\tma H_{m,n}$ may be, and likely is, full rank.
Thus we cannot naively apply Ho and Kalman's theorem to compute the McMillan degree.

\subsection{Lower Bound\label{sec:intro:bound}}
Weyl's theorem provides a way to use Ho and Kalman's theorem
to obtain a lower bound for the McMillan degree.
Recall the rank of any matrix $\ma A \in \C^{m\times n}$
is the number of nonzero singular values;
denoting the $k$th singular value of $\ma A$ as $\sigma_k(\ma A)$
\begin{equation}\label{eq:rank}
	\rank(\ma A) = \sum_{k=1}^{\min (m,n)}  \mathbb{I} [ \sigma_k(\ma A)],
	\qquad 
	\mathbb{I}[\alpha] \coloneqq 
	\begin{cases}
		0, & \alpha <  0 ; \\
		1, & \alpha \ge 0;
	\end{cases}
\end{equation}
where $\mathbb{I}$ is the indicator function.
Using Weyl's theorem \cite[Cor.~7.3.8]{HJ85} provides a bound
connecting the singular values of $\ma H_{m,n}$ and $\tma H_{m,n}$:
\begin{equation}
	|\sigma_k(\tma H_{m,n}) - \sigma_k(\ma H_{m,n})| \le \| \tma H_{m,n} - \ma H_{m,n}\|_2.
\end{equation}
This provides a lower bound on the singular values of $\ma H_{m,n}$:
\begin{equation}
	\sigma_k(\ma H_{m,n}) \ge \sigma_k(\tma H_{m,n}) - \|\tma H_{m,n} - \ma H_{m,n}\|_2.
\end{equation}
Combining this result, equation~\cref{eq:rank}, and \cref{thm:HK}
provides a lower bound on the McMillan degree
\begin{equation}\label{eq:lower_bound}
	\op M(\ve h) \ge \rank(\ma H_{m,n})
		 \ge \sum_{k=1}^{\min(m,n)} \mathbb{I}[ \sigma_k(\tma H_{m,n})  - \|\tma H_{m,n} - \ma H_{m,n}\|_2].
\end{equation}
Although this bound requires computing a quantity we cannot measure---namely the threshold $\|\tma H_{m,n} - \ma H_{m,n}\|_2$---if this threshold is sufficiently small, 
we can visually identify an appropriate approximate threshold as illustrated in \cref{fig:example};
see, e.g., \cite[subsec.~3.5]{Sch00}.
This is necessary even with exact data in $\ma H_{m,n}$ (rounding to floating point accuracy)
as computing the SVD using standard, backward stable algorithms
implies we recover the singular values of a nearby $\tma H_{m,n}$, not those of $\ma H_{m,n}$.

\begin{figure*}
\centering
\begin{tikzpicture}
	\pgfplotstableread{fig_example.dat}\data
	\begin{groupplot}[
		group style = {
			group size = 3 by 1,
			horizontal sep = 0.5em,
			},
		width = 5.3cm,
		height = 5.3cm, 
		ymode = log,
		ymin = 1e-16,
		ymax = 1e2,
		ytickten = {-16, -14, ..., 2},
		xtick = {1, 10, 20, ..., 60},
		xlabel = $k$,
		xmin = 0,
		legend style = {draw = none}, 
	]
		
		\nextgroupplot[title = {$\ve g_{32} = 0$}, ylabel = singular value,]
		\addplot[black, mark=o, only marks] table [x = k, y = sigma_true] {\data};
		\addplot[black, mark=*, only marks, mark size = 1pt] table [x = k, y = sigma_0e+00] {\data};

		\nextgroupplot[title = {$\ve g_{32}\sim \mathcal{N}(\ve 0, 10^{-8} \ma I)$}, yticklabels = {},
			legend to name = fig_example_legend,
			legend entries = {$\sigma_k(\ma H_{17,16})\qquad $,
								$\sigma_k(\tma H_{17,16})\qquad$,
								Theorem~\ref{thm:main} bound on $\|\ma G_{17,16}\|_2$ },
			legend columns = 3,
			legend style={draw=black}
		]
		\addplot[black, mark=o, only marks] table [x = k, y = sigma_true] {\data};
		\addplot[black, mark=*, only marks, mark size = 1pt] table [x = k, y = sigma_1e-04] {\data};
		\addplot[mark=none, colorbrewerA1, thick] coordinates {(0,29.58e-4) (32,29.58e-4)};

		\nextgroupplot[title = {$\ve g_{32}\sim \mathcal{N}(\ve 0, 10^{-2} \ma I)$}, yticklabels = {}]
		\addplot[black, mark=o, only marks] table [x = k, y = sigma_true] {\data};
		\addplot[black, mark=*, only marks, mark size = 1pt] table [x = k, y = sigma_1e-01] {\data};
		\addplot[mark=none, colorbrewerA1, thick] coordinates {(0,29.58e-1) (32,29.58e-1)};

	\end{groupplot}
	\node at (group c2r1.south) [anchor=north, yshift = -2.5em] {\ref{fig_example_legend}};
\end{tikzpicture}

\caption{	An example estimating the McMillan degree using the singular values of $\tma H_{17,16}$.
	Here the true impulse response $\ve h$ was generated from a real state-space system
	with McMillan degree six with
	$\ma A = \diag[0.9 - 0.4i, 0.9 + 0.4i, 0.9 + 0.2i, 0.9 - 0.2i, 0.7, 0.60 ]$,
	$\ve b = \ve 1$, and $\ve c = \ve 1$.
	On the left, even with no noise
	all the singular values of $\tma H_{17,16}$ are nonzero 
	as a result of finite precision computation,
	but it is easy to identify a threshold for computing the rank of $\tma H_{17,16}$.
	In the middle with a moderate amount of noise,
	the magnitude of the trailing singular values has increased,
	but visually we can still identify a threshold for identfying the McMillan degree.
	Note that the bound given in \cref{thm:main} matches visual intuition.
	On the right with a significant amount of noise
	the bound underestimates the McMillan degree as four.
	This is as expected as our result only provides a \emph{lower bound} on the McMillan degree.
	}
\label{fig:example}
\end{figure*}
 
\subsection{Bounding Noise\label{sec:intro:noise}}
In order to make the lower bound on the McMillan degree in~\cref{eq:lower_bound}
practical, we must estimate the threshold $\|\tma H_{m,n} - \ma H_{m,n}\|_2$.
Here we make the assumption that the noise in $\widetilde h_k$ 
is additive and independent \linelabel{line:add_ind} of $h_k$ so that $\widetilde h_k = h_k + g_k$.
Thus the threshold is the $2$-norm of a structured random matrix $\ma G_{m,n}$:
\begin{equation} \label{eq:Gn}
	\ma G_{m,n}
	\coloneqq \begin{bmatrix}
			 {g}_0 &  {g}_1 & \cdots &   {g}_{n-1} \\
			 {g}_1 &  {g}_2 & \cdots & {g}_{n} \\
				\vdots & & & \vdots \\
			 {g}_{m - 1} &  {g}_{m} & \cdots &  {g}_{m+n-2}
		\end{bmatrix}
	= \tma H_{m,n} - \ma H_{m,n}.
\end{equation}
In this paper we construct a probabilistic upper bound on $\|\ma G_{m,n}\|_2$ in \cref{thm:main}
under the assumption that $\lbrace g_k \rbrace_{k=0}^\infty$ samples two variants of Gaussian random noise. 

\subsubsection{Real Gaussian\label{sec:intro:noise:real}} 
The first case considers real-valued Gaussian random noise.
Here we denote the first $N$ entries of $\lbrace g_k \rbrace_{k=0}^\infty$ 
as the vector $\ve g_{N} \in \R^{N}$ and assume $\ve g_{N}$
samples a real-valued multivariate normal distribution with mean zero
and covariance $\ma \Sigma \in \R^{N\times N}$,
denoted $\ve g_{N}\sim \mathcal{N}(\ve 0, \ma \Sigma)$
where $\ma \Sigma$ is a symmetric positive definite (SPD) matrix;
$\ve g_N$ has probability density 
\begin{equation}
	p(\ve g_N) = (2\pi)^{-\frac{N}{2}} 
		(\det \ma \Sigma)^{-\frac12} \exp[-\tfrac12 \ve g_N^\trans \ma \Sigma^{-1} \ve g_N].
\end{equation} 

\subsubsection{Complex Gaussian\label{sec:intro:noise:complex}} 
The second case considers complex-valued Gaussian random noise,
again denoting the first $N$ entries of $\lbrace g_k \rbrace_{k=0}^\infty$
as $\ve g_{N}\in \C^{N}$.
Complex normal distributions require more care to define than their real counterparts.
One approach is to describe $\ve g_N$ in terms of its real and imaginary parts,
\begin{equation}
	\begin{bmatrix}
		\Re \ve g_N \\
		\Im \ve g_N 
	\end{bmatrix}
	\sim \mathcal{N}
	\left(
		\begin{bmatrix} \ve 0 \\ \ve 0 \end{bmatrix},
		\begin{bmatrix}
			\ma \Sigma_{11} & \ma \Sigma_{12} \\
			\ma \Sigma_{12}^\trans & \ma \Sigma_{22}
		\end{bmatrix}
	\right)
	\quad \text{where }
		\begin{bmatrix}
			\ma \Sigma_{11} & \ma \Sigma_{12} \\
			\ma \Sigma_{12}^\trans & \ma \Sigma_{22}
		\end{bmatrix}
	\text{ is SPD}.
\end{equation}
Instead we follow Schreier and Scharf~\cite{SS10}
and characterize $\ve g_N$ via an augmented
complex vector $\uve g_N \in \C^{2N}$ 
containing $\ve g_n$ and its conjugate $\ove g_N$~\cite[sec.~2.1]{SS10}\footnote{
We denote the conjugate of $\ve g$ by $\ove g$
and the complex conjugate transpose of $\ve g$ by $\ve g^*$;
whereas Schreier and Scharf
denote the conjugate of $\ve g$ by $\ve g^*$
and the complex conjugate transpose of $\ve g$ by~$\ve g^{\mathrm H}$.
}
\begin{equation}
\uve g_N \coloneqq \begin{bmatrix} \ve g_N \\ \ove g_N \end{bmatrix}
	= \begin{bmatrix} 
		\ma I & i \ma I \\
		\ma I & -i \ma I
	\end{bmatrix}
	\begin{bmatrix}
		\Re \ve g_N \\
		\Im \ve g_N 
	\end{bmatrix}.
\end{equation}
Now consider the covariance of $\uve g_N$.
As the expected value of $\uve g_N$ is zero, the covariance matrix is simply the expected value of the outer product
$\uve g_n\uve g_n^*$ (cf.~\cite[sec.~2.2]{SS10})
\begin{multline}
	\E( \uve g_N \uve g_N^*)
	=
	\E \left(
	\begin{bmatrix}
		\ma I & i \ma I \\
		\ma I & -i \ma I
	\end{bmatrix}
	\begin{bmatrix}
		\Re \ve g_N \\
		\Im \ve g_N 
	\end{bmatrix}
	\begin{bmatrix}
		\Re \ve g_N \\
		\Im \ve g_N 
	\end{bmatrix}^*
	\begin{bmatrix}
		\ma I & i \ma I \\
		\ma I & -i \ma I
	\end{bmatrix}^*
	\right) 
	 \\
	=\!
	\setlength\arraycolsep{4pt}
	\begin{bmatrix}
		\ma \Sigma_{11} + \ma \Sigma_{22} - i \ma \Sigma_{12} + i \ma \Sigma_{12}^\trans &
		\ma \Sigma_{11} - \ma \Sigma_{22} + i \ma \Sigma_{12} + i \ma \Sigma_{12}^\trans \\
		\ma \Sigma_{11} - \ma \Sigma_{22} - i \ma \Sigma_{12} - i \ma \Sigma_{12}^\trans & 
		\ma \Sigma_{11} + \ma \Sigma_{22} + i \ma \Sigma_{12} - i \ma \Sigma_{12}^\trans 
	\end{bmatrix}
	\!=\! 
	\setlength\arraycolsep{2pt}
	\begin{bmatrix}
		\ma \Gamma & \tma \Gamma \\
		\tma \Gamma^* & \oma \Gamma
	\end{bmatrix}.\!\!\!\! 
\end{multline}
In contrast with real normal distributions, which are completely described by their mean and covariance,
describing $\ve g_N$ requires the mean, 
the Hermitian \emph{covariance matrix} $\ma \Gamma \in \C^{N\times N}$,
and the symmetric \emph{complementary covariance matrix} $\tma \Gamma \in \C^{N\times N}$.

There is a case where specifying a complex normal distribution simplifies
and the resulting random variable acts similar to the real case.
For a generic complex random variable $\ve z$ we say:
\begin{itemize}
\item $\ve z$ is \emph{proper} if the complementary covariance $\E[(\ve z - \E[\ve z])(\ve z - \E[\ve z])^\trans]$ is zero~\cite[Def.~2.1]{SS10};
\item $\ve z$ is \emph{circular} if the probability density of
	 any complex rotation $e^{i\theta}\ve z$ for $\theta \in [0,2\pi)$
	is identical to that of $\ve z$~\cite[Def.~2.4]{SS10} (this requires $\E[\ve z]=\ve 0$). 
\end{itemize}
For a complex normally distributed random variable $\ve z$ with zero mean,
$\ve z$ is proper if and only if $\ve z$ is circular~\cite[Res.~2.11]{SS10}.
Hence circular Gaussian random variables
are completely specified by their covariance $\ma \Gamma \in \C^{N\times N}$.
Here we exclusively consider circular Gaussian random variables,
denoted $\ve g_N\sim \mathcal{CN}(\ve 0, \ma \Gamma)$
where $\ma \Gamma\in \C^{N\times N}$ is Hermitian positive definite;
then $\ve g_N \in \C^N$ has probability density~\cite[Res.~2.5]{SS10}
\begin{equation}
	p(\ve g_N) = \pi^{-n} \det(\ma \Gamma_N)^{-1} \exp[-\ve g_N^*\ma \Gamma^{-1} \ve g_N].
\end{equation}

\subsection{Related Problems\label{sec:intro:related}}
Estimating the McMillan degree via this Hankel matrix approach
is closely related to many problems in system identification and signal processing~\cite{CL15}.
For example, given a complex sinusoidal signal
\begin{equation}
	y(t) = \sum_{k=1}^q \alpha_k e^{\omega_k t} \quad \alpha_k, \omega_k \in \C,
\end{equation}
we can compute the number of components $q$
by considering the McMillan degree of the sequence $\lbrace y(\delta j) \rbrace_{j=0}^\infty$
for some time step $\delta>0$; cf.~\cite[subsec.~2.3]{CL15}.

\subsection{Contributions\label{sec:intro:contributions}}
Here we develop a new probabilistic
upper bound on the 2-norm of a random Hankel matrix $\ma G_{m,n}$ in \cref{thm:main} 
based on a circulant embedding.
Unlike existing results for random Hankel matrices (summarized in \cref{sec:background:structured})
we are able to obtain an upper bound with a fixed probability 
with no unknown constants matching existing asymptotic rate results.
Combined with~\cref{eq:lower_bound},
this upper bound on $\|\ma G_{m,n}\|_2$ allows us to
obtain a lower bound on the McMillan degree of noisy measurements
of the impulse response $\tve h$,
extending Ho and Kalman's result for noisy data.
As illustrated in \cref{sec:examples},
this bound provides a practical estimate of the McMillan degree.
Replacing this probabilistic upper bound on $\|\ma G_{m,n}\|_2$
with an empirical estimate as described in \cref{sec:empirical} provides an even sharper estimate.
Finally, estimating the McMillan degree based on the singular values of a Hankel matrix
compares favorably to \emph{model selection} approaches such as the Akaike Information Criteria (AIC).
Model selection requires 
identifying a minimal realization for each potential McMillan degree,
a process that is both expensive
and prone to identify an unrepresentative local minimum far from the global minimizer.
By using our Hankel matrix approach for estimating the McMillan degree
we avoid this expense and complication.

 \section{Background\label{sec:background}}
Estimating the McMillan degree touches on four distinct domains:
fast Hankel-vector products,
structured random matrices, 
heuristics from engineering practice,
and model order selection.
In the following, we briefly review relevant results from each domain.

\subsection{Fast Hankel Matrix-vector Products\label{sec:background:fast}}
Although $\tma H_{m,n} \in \C^{m\times n}$ is a dense matrix, 
we can exploit the Hankel structure to provide fast matrix-vector products~\cite[sec.~3.4]{Ng04}
and hence accelerate the computation of the SVD.
One approach for fast Hankel vector products is to recognize
a Hankel matrix can be embedded inside a circulant matrix,
which in turn can be diagonalized by the discrete Fourier transform (DFT) matrix.
This allows the product $\tma H_{m,n} \ve x$ to be computed 
using only $\order(N \log N)$ operations where $N=m+n-1$,
rather than the $\order(mn)$ operations normally required. 
These inexpensive inner products can then accelerate the computation of the SVD 
when using an iterative eigensolver like ARPACK~\cite{LSY98},
with the leading $k$ singular values be computed in approximately $\order(k N \log N)$ operations.

\subsection{Structured Random Matrices\label{sec:background:structured}}
The spectral properties of structured random matrices have only started to be explored in the past two decades.
The distribution of the singular values of a random Hankel matrix (and hence the 2-norm) 
was posed as an open problem in a 1999 paper by Bai~\cite{Bai99}.
Byrc, Dembo, and Jiang were the first to establish the limiting spectral distribution for Hankel matrices
with independent and identically distributed (iid) Gaussian entries in 2006~\cite{BDJ06}.
The next year, Meckes provided bounds on the distribution of the 2-norm
under weaker assumptions that entries are uniformly subgaussian, independent, but not necessarily identically distributed~\cite{Mec07}.
Combining Meckes' Theorem~1 and Theorem~3
we know the growth rate of $\E \|\ma G_{n,n}\|_2$ as a function of $n$;
assuming the entries of $\ma G_{n,n}$ are iid Gaussian random variables
with zero mean and unit variance,
then there exists $0<c_1< c_2$ such that  
\begin{equation}\label{eq:asymptotic_growth}
		c_1 \sqrt{n\log n} \le \E \|\ma G_{n,n}\|_2  \le c_2\sqrt{n \log n} \quad \forall n >0.
\end{equation}
Similar results were established under even weaker constraints for the distribution of the entries
by Adamczak~\cite{Ada10} and Nekrutkin~\cite{Nek13}; the latter also treated non-square Hankel matrices.
Note that although our results require a more restrictive assumption that entries of $\ma G_{n,n}$
sample a multivariate Gaussian distribution,
we provide a different result: a computable probabilistic upper bound on $\|\ma G_{n,n}\|_2$.

\subsection{Heuristics for Estimating McMillan Degree}
Although rigorous estimates of the 2-norm of 
a random Hankel matrix have only been available for the past two decades,
many authors in the 1970s, 1980s, 1990s 
recognized that the singular values of $\tma H_{m,n}$ could
be used to infer the McMillan degree.
For example, in 1985 Juang and Pappa suggested picking a threshold manually
to separate singular values into those associated with $\ma H_{m,n}$
and those associated with noise~\cite[p.622]{JP85}---a process that as illustrated in \cref{fig:example} sometimes yields 
an obvious choice, but that sometimes can be misleading.
This manual selection approach also appears in more recent work using matrices related to $\ma H_{m,n}$;
see, e.g.,\cite[\S16.3]{Lju99}, \cite{ZZXH05}, and \cite{YL07}.
Other authors have attempted to provide estimates of $\|\ma G_{m,n}\|_2$ 
to select this threshold in~\cref{eq:lower_bound}.
For example, Holt and Antill bounded the norm of a Hankel matrix by its Frobenius norm~\cite[eq.~(19)]{HA77}.
Assuming $\ve g_{2n-1} \sim \mathcal{N}(\ve 0, \epsilon\ma I)$,
\begin{equation}
	\|\ma G_{n,n}\|_2 \le \|\ma G_{n,n}\|_\fro
	\quad \Longrightarrow \quad
	\E [ \|\ma G_{n,n}\|_2]
	\le \E [ \|\ma G_{n,n}\|_\fro]
	=   \sqrt{n^2\E[ g_0^2]} = n\epsilon
\end{equation}
However this bound is far too conservative:
from~\cref{eq:asymptotic_growth} we know $\|\ma G_{n,n}\|_2$
grows with $n$ like $\order(\sqrt{n \log n})$,
whereas this bound is $\order(n)$.
Another threshold that has been suggested
when $\ve g_{2n-1} \sim \mathcal{N}(\ve 0, \epsilon \ma I)$
is $\epsilon \sqrt{n}$;
see, e.g., \cite[eq.~(4.3)]{GM87} and \cite[\S IV.C]{VDS93}.
This is based on the expected value of $\tma H_{n,n}^*\tma H_{n,n}$
\begin{equation}
	\begin{split}
	\E[\tma H^*_{n,n}\tma H_{n,n}] 		&= \ma H_{n,n}^*\ma H_{n,n} + \E[ \ma G_{n,n}^*\ma H_{n,n}] 
			+ \E[ \ma H_{n,n}^*\ma G_{n,n}] + \E[\ma G_{n,n}^*\ma G_{n,n}] \\
		&= \ma H_{n,n}^*\ma H_{n,n} + \epsilon^2 n \ma I
	\end{split}
\end{equation}
whose eigenvalues are all shifted upwards by $\epsilon^2n$;
hence the singular values of the matrix square root of $\E[\tma H_{n,n}^*\tma H_{n,n}]$
are shifted upwards by $\epsilon\sqrt{n}$.
However this threshold makes a mistake interchanging expectation and the eigenvalues:
the eigenvalues of $\E[\tma H_{n,n}^*\tma H_{n,n}]$ \linelabel{line:switch}
are \emph{not} the expected eigenvalues of $\tma H_{n,n}^*\tma H_{n,n}$.
The result is a threshold that is too permissive as it grows like $\order(\sqrt{n})$
whereas we should expect $\order(\sqrt{n\log n})$.

\subsection{Model Selection\label{sec:background:model}}
Model selection provides an alternative perspective on estimating the McMillan degree
using generic statistical tools for selecting the most parsimonious model among a set of candidate models.
In the context of estimating the McMillan degree,
the candidate models are realizations consisting of matrices
$\ma A \in \C^{q\times q}$ and vectors $\ve b, \ve c \in \C^q$
for differing dimensions $q$.
There are a large number of different criteria for selecting the most parsimonious model, see, e.g.,~\cite{BA02}.
Here we focus on information theoretic approaches
which score candidate models on both likelihood and number of parameters.
The Akaike Information Criteria (AIC)~\cite{Aka74} is one such popular model selection criteria 
where the score of each model is proportional to the number of free parameters minus the log-likelihood.
In our context, 
for either real $\ve g_{n} \sim \mathcal{N}(\ve 0, \ma \Sigma_{n})$ or 
complex circular $\ve g_{n}\sim \mathcal{CN}(\ve 0, \ma \Sigma_{n})$
Gaussian random noise,
the AIC score for a model of degree $q$ is:
\begin{equation}
	\text{AIC}(q) \propto 2\min_{\substack{\ma A\in \C^{q\times q} \\ \ve b,\ve c\in \C^q} }
						\left\|\ma \Sigma^{-\frac12} 
						\left(
							\begin{bsmallmatrix}
								\widetilde h_0 \\ \widetilde h_1 \\ \vdots \\ \widetilde h_{n-1}
							\end{bsmallmatrix}
							-
							\begin{bsmallmatrix}
								\ve c^* \ma A^0 \ve b \\ \ve c^*\ma A^1 \ve b \\ \vdots \\ \ve c^*\ma A^{n-1} \ve b
							\end{bsmallmatrix}
						\right)
						\right\|_2^2 + 4q
				+ \text{constant}.
\end{equation}
The second term in the AIC encodes the number of real degrees of freedom in the model.
Although $\ma A$, $\ve b$, and $\ve c$ have a collective $q^2+2q$ degrees of freedom,
there are only effectively $4q$ degrees of freedom.
Without loss of generality we can assume $\ve c=\ve 1$
and that $\ma A$ is diagonal as non-diagonalizable matrices 
are nowhere dense in $\C^{q\times q}$~\cite[p.~2739]{DGA10};
this leaves $2q$ complex parameters or $4q$ real parameters.
Then the AIC selects the $q$ minimizing $\text{AIC}(q)$.
The challenge with this approach is its expense:
for each candidate McMillan degree a minimal realization $\lbrace \ma A, \ve b, \ve c\rbrace$
must be constructed.

 \section{Random Hankel Matrix 2-Norm Bound\label{sec:bound}}
We now establish our main result:
a probabilistic upper bound on the 2-norm of a random Hankel matrix
whose entries are drawn from a multivariate normal distribution.

\begin{theorem}\label{thm:main}
	Suppose $\ve g_N\in \C^N$ is a random variable
	and $\ma G_{m,n}\in \C^{m \times n}$ is a Hankel matrix
	constructed from $\ve g_N$ as in~\cref{eq:Gn} where $N=m+n-1$,
	then
	\begin{equation}\label{eq:Gnorm}
		\| \ma G_{m,n}\|_2 \le \alpha \sqrt{N} \text{ with probability $p(\alpha)$}
	\end{equation}
	where $p(\alpha)$ depends on the distribution of $\ve g_N$:
	\begin{subequations}\label{eq:probs}
	\begin{align}
		\label{eq:probs:real}
		&\text{if $\ve g_N \sim \mathcal{N}(\ve 0, \ma I)$, then}
			& 
				p(\alpha) &= 
				\begin{cases}
			\erf(\alpha/2)^{\phantom{2}}(1-e^{-\alpha^2/2})^{(N-1)/2}, & N \text{ odd};\\
			\erf(\alpha/2)^2(1-e^{-\alpha^2/2})^{N/2-1}, & N \text{ even};
				\end{cases} \\
		\label{eq:probs:complex}
		&\text{if $\ve g_N \sim \mathcal{CN}(\ve 0, \ma I)$, then}
				& p(\alpha) &= (1- e^{-\alpha^2/2})^N; \\
		&\text{if $\ve g_N\sim\mathcal{N}(\ve 0, \ma \Sigma)$, then}
			&	p(\alpha) &= \gamma(N/2, \alpha^2/(2\| \ma \Sigma^{\frac12}\|_2^2))/\Gamma(N/2); \\
		&\text{if $\ve g_N\sim \mathcal{CN}(\ve 0, \ma \Sigma)$, then}
			&	p(\alpha) &= \gamma(N, \alpha^2/\| \ma \Sigma^{\frac12}\|_2^2 )/\Gamma(N)
	\end{align}
	\end{subequations}
	where $\Gamma$ denotes the Gamma function, $\Gamma(s) \coloneqq \int_0^\infty t^{s-1} e^{-t} \D t$,
	$\gamma$ is the lower incomplete gamma function, $\gamma(s,x) \coloneqq \int_0^x t^{s-1}e^{-t} \D t$,
	and $\erf$ is the error function, $\erf(x) \coloneqq 2\pi^{-1/2}\int_0^x e^{-t^2} \D t$.
\end{theorem}

We are able to state this result for any rectangular Hankel matrix $\ma G_{m,n}$
as the first component of the proof---a circulant embedding to obtain a bound in terms of the DFT of $\ve g_N$---yields the same bound for any Hankel matrix the same generating data of length $N = m+n-1$.
The second component then takes this bound and generates a probabilistic upper bound 
assuming a particular distribution for $\ve g_N$.

\subsection{Asymptotic Growth\label{sec:bound:growth}}
Before proving this result, we ask:
does $\|\ma G_{n,n}\|_2$ grow at the same rate as $n\to\infty$
as the bound provided by Meckes~\cite[Thm.~3]{Mec07}, namely $\order(n\log n)$?
This is true for the circular complex normal case~\cref{eq:probs:complex}.
For a fixed probability $\tau\in (0,1)$,
then the $\alpha$ satisfying $\tau = p(\alpha)$ is
\begin{equation}\label{eq:growth_bnd}
	\alpha = \sqrt{-2\log(1-\tau^{\frac{1}{N}})}
	= \sqrt{2 \log N  - \log (\log \tau)^2 + \order(N^{-1})} = 
	\order(\log N)
\end{equation}
as $N\to \infty$. 
Here we used a Taylor expansion of the exponential in $\tau^{\frac{1}{N}} = \exp[\log[\tau^{\frac{1}{N}}]]$
to obtain this estimate.
Hence in~\cref{eq:Gnorm}, $\|\ma G_{n,n}\|_2 = \order(\sqrt{N \log N}) = \order(n\log n)$
with probability $\tau$ when $\ve g_N \sim \mathcal{CN}(\ve 0, \ma I)$.

\begin{figure*}
\centering

\begin{tikzpicture}
\begin{groupplot}[
	group style = {group size = 2 by 1, horizontal sep = 4em},
	width = 0.51\linewidth,
	height = 0.46\linewidth,
	xlabel = $N$,
	ylabel = {$\beta$ s.t. $\prob [ \|\ma G_{m,n}\|_2 < \beta)] = \tau$},
	xmode = log, ymode = log,
	xmin = 0.9, xmax = 1.2e6,
	ymin = 0.2, ymax = 1e4,
	ytickten = {0,1, 2, 3, 4}
	]

	\nextgroupplot[ title = {$\ve g_N \sim \mathcal{N}(\ve 0, \ma I)$}]
	\addplot[black,thin] table [x=n, y = bnd99] {fig_random_hankel_real.dat} 
		[anchor = east] node [pos=0.0, rotate =0, xshift = 5pt] {\scriptsize $99\%$}
		[anchor = south] node [pos = 0.6, rotate = 31] {\scriptsize \cref{thm:main}a};
	\addplot[black, thick] table [x=n, y = bnd50] {fig_random_hankel_real.dat}
		[anchor = east] node [pos=0.0, rotate =0, xshift = 5pt] {\scriptsize $50\%$};
	\addplot[black,thin] table [x=n, y = bnd1] {fig_random_hankel_real.dat}
		[anchor = east] node [pos=0.0, rotate =0, xshift = 5pt] {\scriptsize $1\%$};
	\addplot[colorbrewerA1, thick, only marks, mark size =3, mark=-, error bars/.cd, y dir = plus, y explicit]
		 table [x=n, y = p50, y error expr = \thisrow{p99} - \thisrow{p50} ] {fig_random_hankel_real.dat};
	\addplot[colorbrewerA1, thick, only marks, mark=-, mark size =3, error bars/.cd, y dir = minus, y explicit]
		 table [x=n, y = p50, y error expr = \thisrow{p50} - \thisrow{p1} ] {fig_random_hankel_real.dat};
	\addplot[colorbrewerA1, only marks, mark=*, mark size=1] table[x=n, y=p50] {fig_random_hankel_real.dat};
	
	\addplot[colorbrewerA2, thin, densely dashed] table [x=n, y = p1r] {fig_random_hankel_l2.dat};
	\addplot[colorbrewerA2, thick, densely dashed] table [x=n, y = p50r] {fig_random_hankel_l2.dat};
	\addplot[colorbrewerA2, thin, densely dashed] table [x=n, y = p99r] {fig_random_hankel_l2.dat}
		[anchor = south] node [pos = 0.3, rotate = 47] {\scriptsize \cref{thm:main}c};

	\addplot[densely dotted, black, domain = 4:1e6] {1e-1*sqrt(x * ln(x))}
		[anchor = north] node [pos = 0.5, rotate =33] {\scriptsize $\sqrt{N \log N}$}; 	
	
	\nextgroupplot[ title = {$\ve g_N \sim \mathcal{CN}(\ve 0, \ma I)$}]
	\addplot[black,thin] table [x=n, y = bnd99] {fig_random_hankel_complex.dat} 
		[anchor = east] node [pos=0.0, rotate =0, xshift = 5pt] {\scriptsize $99\%$}
		[anchor = south] node [pos = 0.6, rotate = 31] {\scriptsize \cref{thm:main}b};
	\addplot[black, thick] table [x=n, y = bnd50] {fig_random_hankel_complex.dat}
		[anchor = east] node [pos=0.0, rotate =0, xshift = 5pt] {\scriptsize $50\%$};
	\addplot[black,thin] table [x=n, y = bnd1] {fig_random_hankel_complex.dat}
		[anchor = east] node [pos=0.0, rotate =0, xshift = 5pt] {\scriptsize $1\%$};

	\addplot[colorbrewerA2, thin, densely dashed] table [x=n, y = p1c] {fig_random_hankel_l2.dat};
	\addplot[colorbrewerA2, thick, densely dashed] table [x=n, y = p50c] {fig_random_hankel_l2.dat};
	\addplot[colorbrewerA2, thin, densely dashed] table [x=n, y = p99c] {fig_random_hankel_l2.dat}
		[anchor = south] node [pos = 0.3, rotate = 47] {\scriptsize \cref{thm:main}d};

	\addplot[colorbrewerA1, thick, only marks, mark size =3, mark=-, error bars/.cd, y dir = plus, y explicit]
		 table [x=n, y = p50, y error expr = \thisrow{p99} - \thisrow{p50} ] {fig_random_hankel_complex.dat};
	\addplot[colorbrewerA1, thick, only marks, mark=-, mark size =3, error bars/.cd, y dir = minus, y explicit]
		 table [x=n, y = p50, y error expr = \thisrow{p50} - \thisrow{p1} ] {fig_random_hankel_complex.dat};
	\addplot[colorbrewerA1, only marks, mark=*, mark size=1] table[x=n, y=p50] {fig_random_hankel_complex.dat};
	
	\addplot[densely dotted, black, domain = 4:1e6] {1e-1*sqrt(x * ln(x))}
		[anchor = north] node [pos = 0.5, rotate =31] {\scriptsize $\sqrt{N\log N}$}; 	
\end{groupplot}
\end{tikzpicture}
\caption{	A comparison of the upper bounds from \cref{thm:main} 
	for $\|\ma G_{m,n}\|_2$ for increasing $N$ where $m = \lceil \frac{N-1}{2}\rceil$
	and $n = \lfloor \frac{N-1}{2}\rfloor$.
	The solid curves show the bounds from \cref{thm:main}
	which hold when the covariance matrix of $\ve g_N$ is the identity matrix;
	the dashed lines show bound allowing other covariance matrices.
	The curves show the 1st, 50th, and 99th percentiles.
	The red bars show an empirical estimate of $\|\ma G_{m,n}\|_2$ based on $10^3$ Monte Carlo samples,
	showing the bound $\beta$ that holds with probability $\tau$;
	the bars similarly correspond to the 1st, 50th, and 99th percentiles.
}
\label{fig:random_hankel}
\end{figure*}
 
\Cref{fig:random_hankel} compares an empirical estimate the distribution of $\|\ma G_{n,n}\|_2$
to the bounds provided by \cref{thm:main}.
We observe that both the real and circular complex normal distribution bounds
in \cref{eq:probs:real} and \cref{eq:probs:complex}
match the expected asymptotic growth rate of $\order(\sqrt{n\log n})$.
Moreover, for these two cases, the bound is only approximately $2.5$ times larger than the empirical estimate.

\subsection{Circulant Embedding Bound}
The first step in establishing Theorem~\ref{thm:main}
bounds $\|\ma G_{m,n}\|_2$ by embedding $\ma G_{m,n}$.
This circulant matrix is diagonalized by 
the discrete Fourier transform matrix (DFT) allowing us to obtain its 2-norm.
Although this circulant embedding technique has long been used for fast Hankel matrix-vector products~\cite[sec.~3.4]{Ng04},
this is, to the best of our knowledge,
the first time this embedding technique has been used to obtain bounds 
on the norm of a Hankel matrix.

\begin{lemma}\label{lem:embedding}
	Suppose $\ve g_N$ and $\ma G_{m,n}$ are defined as in~\cref{thm:main},
	then
	\begin{equation}
		\| \ma G_{m,n} \|_2 \le \sqrt{N} \| \ma F_{N} \ve g_N\|_\infty
	\end{equation}
	where $[\ma F_N]_{j,k} = N^{-\frac12} e^{-2\pi i jk/N}$ is the DFT matrix.
\end{lemma}
\begin{proof}
	Let $\ma C_N \in \C^{N\times N}$ be a circulant matrix~\cite[\S0.9.6]{HJ85} whose first column is $\ve g_N$
	and recalling $N = m + n -1$,
	\begin{equation}\label{eq:circ}
		\ma C_N =  \begin{bmatrix}
			g_0 & g_{m+n-2} & \ldots & g_{n-1} & \! g_{n-2} & \ldots & g_1  \\
			g_1 & g_0 & \ldots & g_{n} & g_{n-1} & \ldots & g_2  \\
			\vdots & & \ddots & \vdots & \vdots & & \vdots \\
			\ \tikzmark{left}\ \ \  g_{m-1}\ \ \  & g_{m-2} & \ldots & g_0 & g_{m+n-2} & \ldots & g_m \\
			g_{m} & g_{m-1} & \ldots & g_1 & g_0 & \ldots  &  g_{m+1}\\
			\vdots & & & \vdots & \vdots & \ddots & \vdots \\
			g_{m+n-2} & g_{m+n-3} & \ldots & g_{n-1}\ \tikzmark{right} & \! g_{n-2} & \ldots & g_0
		\end{bmatrix}
	\DrawBox[thick]
	\end{equation}
	Note the Hankel matrix $\ma G_{m,n}\in \C^{m\times n}$ appears in the boxed region of $\ma C_N$ with reversed columns.
	Hence the multiplication $\ma G_{m,n}\ve x_n$ can be written as
	\begin{equation}\label{eq:embedding}
		\ma G_{m,n} \ve x_n =
			\begin{bmatrix}
				\ma 0  &
				 \ma I_{m}
			\end{bmatrix}
			\ma C_N
			\begin{bmatrix}
			\ma J_n  \\ \ma 0
			\end{bmatrix}
			\ve x_n
	\end{equation}
	where $\ma I_m \in \C^{m\times m}$ is the identity matrix
	and $\ma J_m \in \C^{n \times n}$ is the identity matrix with columns reversed.
	Then, as the matrix 2-norm is induced by the vector 2-norm,
	\begin{align}
		\| \ma G_{m,n} \|_2 &\coloneqq 
			\max_{\ve x_n\in \C^n\backslash \lbrace 0 \rbrace} 
				\frac{ \|\ma G_{m,n} \ve x_n\|_2}{\| \ve x_n\|_2} 
		= \max_{\ve x_n\in \C^n\backslash \lbrace 0 \rbrace} 
			\frac{\left\|
			\begin{bmatrix}
				\ma 0  & \ma I_m
			\end{bmatrix}
			\ma C_N
			\begin{bmatrix}
			\ma J_n  \\ \ma 0
			\end{bmatrix}
			\ve x_n
			\right\|_2}{\|\ve x_n\|_2}\\
		&\le \max_{\ve x_n\in \C^n\backslash \lbrace 0 \rbrace} 
			\frac{\left\|
			\ma C_N
			\begin{bmatrix}
			\ma J_n  \\ \ma 0
			\end{bmatrix}
			\ve x_n
			\right\|_2}{\|\ve x_n\|_2}
		\le \max_{\ve y_N \in \C^N\backslash \lbrace 0 \rbrace }
			\frac{\left\| \ma C_N\ve y_N\right\|_2}{\| \ve y_N\|_2}
		= \| \ma C_N\|_2.
	\end{align}
	Finally, to bound the norm of $\ma C_N$
	we note that since  $\ma C_N$ is a circulant matrix, it has spectral decomposition~\cite[eq.~(3.27)]{Ng04},
	\begin{equation}
		\ma C_N =  \ma F_N^* \ma \Lambda_N \ma F_N, \qquad \ma \Lambda_N  = \sqrt{N} \diag(\ma F_N \ve g_N),
	\end{equation}
	and then, as the 2-norm is unitarily invariant,
	\begin{equation}
		\|\ma C_N \|_2 = \|\ma F_N^*\ma \Lambda_N \ma F_N\|_2 = \| \ma \Lambda_N \|_2
		= \sqrt{N}\| \ma F_N \ve g_N\|_\infty.
	\end{equation}
\end{proof}
 \subsection{Bounds on Noise\label{sec:bound:lemma}}
We now seek to bound $\|\ma F_N \ve g_N\|_\infty$
for four different distributions associated with $\ve g_N$,
corresponding to the four cases in \cref{eq:probs}. 
Here we denote the probability of an expression being true by $\prob$;
e.g., the probability of $z \le \alpha$ being true for some random variable $z$ is
$\prob[ z \le \alpha] \coloneqq \E[\mathbb{I}[\alpha - z]]$.

\begin{lemma}\label{lem:complex}
	Suppose $\ve g_N\sim \mathcal{CN}(\ve 0, \ma I)$ and $\alpha \ge 0$,
	then
	\begin{equation}
		\prob [ \, \| \ma F_N \ve g_N\|_\infty \le \alpha ] = (1 - e^{-\alpha^2/2})^{N}.
	\end{equation}
\end{lemma}
\begin{proof}
	We begin by characterizing $\ma F_N \ve g_N$.
	Note $\E[\ma F_N \ve g_N]=\ve 0$ and 
	hence the covariance and complementary covariance matrices are
	\begin{equation}
		\E[\ma F_N \ve g_N \ve g_N^*\ma F_N^*] = \ma F_N\ma I \ma F_N^* = \ma I
		\qquad 
		\E[\ma F_N \ve g_N \ve g_N^\trans\ma F_N^\trans] = \ma F_N\ma 0 \ma F_N^\trans = \ma 0
	\end{equation} 
	where the second statement follows as $\ve g_N$ is circular.
	Hence $\ma F_N\ve g_N$ is a circular Gaussian random variable 
	with $\ma F_N \ve g_N \sim \mathcal{CN}(\ma 0, \ma I)$; cf.~\cite[subsec.~2.3.1]{SS10}.
	As such, the $k$ entry of $\ma F_N \ve g_N$ is independent of the $\ell$th entry
	when $k \ne \ell$
	and hence 
	\begin{multline}\label{eq:prob_chain}
		\prob[ \, \|\ma F_N\ve g_N]\|_\infty \le \alpha]
		 =\prob[ \, \max_k |\ve e_k^*\ma F_N\ve g_N | \le \alpha]  
		= \prod_{k=0}^{N-1} \prob[ \, |\ve e_k^* \ma F_N \ve g_N | \le \alpha]
	\end{multline}
	where $\ve e_k$ denotes the $k$th column of the identity.
	Note $\ve e_k^*\ma F_N \ve g_N$ has the distribution
	\begin{equation}
		\ve e_k^*\ma F_N \ve g_N
			\sim \mathcal{CN}(0, 1 ).
	\end{equation}
	Hence $| \ve e_k^* \ma F_N \ve g_N |$ follows a Rayleigh distribution \linelabel{line:rayleigh}
	(i.e., $\chi_2$, a $\chi$-distribution with two degrees of freedom)
	with the cumulative density function~\cite[eq.~(2.74)]{SS10}
	\begin{equation}
		\prob[ \, | \ve e_k^*\ma F_N \ve g_N | \le \alpha] = 1 - e^{-\alpha^2/2}.
	\end{equation}
	Combining this with~\cref{eq:prob_chain} provides the desired probability.
\end{proof}

The analogous result for real Gaussian random variables $\ve g_N \sim \mathcal{N}(\ve 0, \ma I)$
requires additional care as the entries of $\ma F_N \ve g_N$ are no longer independent---half of the entries are conjugates of the other half.

\begin{lemma}\label{lem:real}
	Suppose $\ve g_N \sim \mathcal{N}(\ve 0, \ma I)$
	and $\alpha \ge 0$, then
	\begin{equation}
		\prob[\, \| \ma F_N \ve g_N\|_\infty \le \alpha]
		=
		\begin{cases}
			\erf(\alpha/2)^{\phantom{2}}(1-e^{-\alpha^2/2})^{(N-1)/2}, & N \text{ odd};\\
			\erf(\alpha/2)^2(1-e^{-\alpha^2/2})^{N/2-1}, & N \text{ even}.
		\end{cases}
	\end{equation}
\end{lemma}
\begin{proof}
	First, we write the real random variable $\ve g_N $
	as a function of the complex circular normal variable $\ve z_N \sim \mathcal{CN}(\ve 0, \ma I)$:
	\begin{equation}
		\ve g_N = 2^{-\frac12}( \ve z_N + \conj{\ve z_N}) \sim \mathcal{N}(\ve 0, \ma I).
	\end{equation}
	Then, defining $\ve w_N := \ma F_N \ve z_N$, 
	\begin{align}
		\ma F_N \ve g_N = 2^{-\frac12}( \ma F_N \ve z_N + \ma F_N \conj{\ve z_N}) 
			= 2^{-\frac12}( \ve w_N + \ma F_N \ma F_N^\trans \conj{\ve w_N}).
	\end{align}
	Above, the matrix $\ma F_N \ma F_N^\trans$ has the form
	\begin{align}
		\ma F_N\ma F_N^\trans = \begin{bmatrix}
			1 & \ve 0^\trans \\
			\ve 0 & \ma J_{N-1}
		\end{bmatrix}
		\in \R^{N\times N},
	\end{align}
	where $\ma J_{N-1}$ is the reversed identity matrix.
	Thus, the entries of $\ma F_N \ve g_N$ are: 
	\begin{equation}\label{eq:Fg}
		\ve e_k^*\ma F_N \ve g_N
		= 2^{-\frac12}  \ve e_k^*\left(
			\begin{bsmallmatrix} w_0 \\ w_1 \\ \vdots \\ w_{N-1} 
			\end{bsmallmatrix}
			+ \begin{bsmallmatrix} \conj{w_0} \\ \conj{w_{N-1}} \\ \vdots \\ \conj{w_1} \end{bsmallmatrix}\right)
		=\begin{cases}
			 2^{\frac12} \Re[ w_0], & k = 0;\\
			2^{-\frac12}( w_k + \conj{w_{N-k}}) & k \ne 0.
		\end{cases}
	\end{equation}
	Then since $\ve w_N = \ma F_N\ve z_N \sim \mathcal{CN}(\ve 0, \ma I)$,
	each entry of $\ma F_N \ve g_N$ is distributed like
	\begin{align}
		\ve e_k^*\ma F_N \ve g_N \sim
		\begin{cases}
			\mathcal{N}(0, 2), & k = 0 \text{ or } k = N/2; \\
			\mathcal{CN}(0,1), & \text{otherwise};
		\end{cases}
	\end{align}
	with cumulative density functions 
	\begin{align}\label{eq:real_entry}
		\prob[ \, |\ve e_k^*\ma F_N \ve g_N| \le \alpha) &= 
		\begin{cases}
			\erf(\alpha/2), & k =0 \text{ or } N/2; \\
			1 - e^{-\alpha^2/2}, & \text{otherwise}.
		\end{cases}
	\end{align}
	Then since the first $\lfloor N/2\rfloor$ entries of $\ma F_N\ve g_N$ are independent
	of each other and the remaining are fully determined by this first half, cf.~\cref{eq:Fg},
	we have 
	\begin{align}
		\prob[\, \| \ma F_N \ve g_N\|_\infty \le \alpha ] =
		\prob\left[\, \max_{k=0,\ldots,\lfloor N/2\rfloor} |\ve e_k^*\ma F_N\ve g_N| \le \alpha\right]
		= \prod_{k=0}^{\lfloor N/2\rfloor}
			\prob [\, |\ve e_k^*\ma F_N\ve g_N| \le \alpha].
	\end{align}
	Then using the entrywise expression~\cref{eq:real_entry}
	we obtain the desired bound.
\end{proof}

Unfortunately we have been unable to find satisfying bounds
on $\|\ma F_N\ve g_N\|_\infty$
when $\ve g_N$ has covariance matrix that is not an identity matrix.
Suppose $\ve g_N \sim \mathcal{CN}(\ve 0, \ma \Sigma)$ for some Hermitian positive definite $\ma \Sigma$;
then $\ma F_N\ve g_N\sim \mathcal{CN}(\ve 0, \ma F_N\ma \Sigma \ma F_N^*)$.
As the entries of $\ma F_N \ve g_N$ are now correlated, 
we cannot separate the probabily of the max into the product of probabilities~\cref{eq:prob_chain}.
The following two bounds provide guideance in this case, 
using the equivalence of finite dimensional norms
and the fact $\ma F_N$ is a unitary matrix:
\begin{equation}\label{eq:sloppy}
	\| \ma F_N \ve g_N\|_\infty \le \|\ma F_N \ve g_N\|_2 = \| \ve g_N\|_2.
\end{equation}
Although this provides a bound,
as evidenced in \cref{fig:random_hankel},
it does not achieve the expected asymptotic growth rate of $\order(\sqrt{N\log N})$
as $N\to \infty$.

\begin{lemma}
	Suppose $\ve g_N \sim \mathcal{N}(\ve 0, \ma \Sigma)$ 
	where $\ma \Sigma\in \R^{N\times N}$ is symmetric positive definite
	and $\alpha \ge 0$, then
	\begin{equation}
		\prob[ \, \|\ma F_N \ve g_N\|_\infty \le \alpha \|\ma \Sigma^{\frac12}\|_2 \, ]
		=1 - \Gamma(N/2)^{-1}\gamma(N/2,\alpha^2/2).
	\end{equation}
\end{lemma}
\begin{proof}
	Writing $\ve g_N = \ma \Sigma^{\frac12}\ve w_N$
	where $\ve w \sim\mathcal{N}(\ve 0, \ma I)$, then
	invoking \cref{eq:sloppy},
	\begin{align}
		\|\ma F_N\ve g_N\|_\infty 
		= \|\ma F_N \ma \Sigma^{\frac12} \ve w_N\|_\infty
		\le \|\ma F_N \ma \Sigma^{\frac12} \ve w_N\|_2
		\le \|\ma \Sigma^{\frac12}\|_2 \|\ve w_N\|_2.
	\end{align}	
	The term $\|\ve w_N\|_2$ samples a $\chi$-distribution with $n$ degrees of freedom
	and the result follows from this density's cumulative distribution.
\end{proof}

The proof for the complex case is identical
except that the $\chi$-distribution has a total of $2N$ degrees of freedom,
with half coming from the real part and half from the imaginary part. 
\begin{lemma}
	Suppose $\ve g_N\sim \mathcal{N}(\ve 0, \ma \Sigma)$
	where $\ma \Sigma \in \C^{N\times N}$ is a Hermitian positive definite matrix
	and $\alpha \ge 0$,
	then
	\begin{equation}
		\prob[ \, \|\ma F_N \ve g_N\|_\infty \le \alpha \|\ma \Sigma^{\frac12}\|_2 \, ]
		 = 1 - \Gamma(N)^{-1}\gamma(N,\alpha^2/2).
	\end{equation}
\end{lemma}

 \section{McMillan Degree Lower Bound\label{sec:order}}
Having provided the bound on the norm of a random Hankel matrix $\|\ma G_{m,n}\|_2$
in \cref{thm:main},
we now formally state the McMillan degree lower bound 
based on this result.

\begin{theorem}\label{thm:order}
	Suppose $\ve h \in \ell_1(\mathbb{N})$
	is the impulse response of a system.
	Given noisy measurements $\widetilde h_k = h_k + g_k$
	constructed into a Hankel matrix $\tma H_{m,n}\in \C^{m\times n}$,
	the McMillan degree of $\ve h$ is bounded below by
	\begin{equation}
		\op M(\ve h) \ge \sum_{k=1}^{\min(m,n)}
			\mathbb{I}[ \, \sigma_k(\tma H_{m,n}) - \alpha \sqrt{m+n-1} \, ]
		\quad \text{with probablity} \quad
		p(\alpha)
	\end{equation}
	where $p(\alpha)$ 
	depends on the distribution of $\ve g_{m+n-1}$
	as given in~\cref{eq:probs}.
\end{theorem}
\begin{proof}
	From \cref{eq:lower_bound} and \cref{eq:Gn} 
	\begin{equation}
		\op M(\ve h) \ge  
			\sum_{k=1}^{\min(m,n)} \mathbb{I}[\, \sigma_k(\tma H_{m,n}) - \|\ma G_{m,n}\|_2] 
	\end{equation}
	From \cref{thm:main} we obtain a probabilistic upper bound on $\|\ma G_{m,n}\|_2$,
	which in turn provides a probabilistic lower bound on $\op M(\ve h)$.
\end{proof}

 \section{Empirical Bound\label{sec:empirical}}
Before continuing to the numerical experiments,
we note that we need not necessarily rely on an 
exact, probabilistic upper bound of $\|\ma G_{m,n}\|_2$.
Instead, as it is inexpensive to compute 
the 2-norm of a Hankel matrix (see \cref{sec:background:fast}),
we can instead sample many realisations of the noise
to estimate the cumulative density function 
associated with the 2-norm of this Hankel matrix.
The advantage of this approach is it provides
sharper estimates of $\|\ma G_{m,n}\|_2$ 
and is applicable to a wider variety of distributions of $g_k$.
However because this is an empirical estimate, 
we cannot provide the guarantees as in \cref{thm:order}.
 \section{Numerical Examples\label{sec:examples}}
Here we provide two examples of our McMillan degree lower bound:
one with complex valued data with a system known McMillan degree
and another with real data with a highly reducible system.
In these examples we compute the AIC score
using HSVD~\cite{BBO87} to estimate the (approximate) optimal model parameters of each candidate McMillan degree.
Following the principles of reproducible research,
code for constructing these examples is available at {\tt \url{https://github.com/jeffrey-hokanson/McMillanDegree}}.

\subsection{Complex Valued Data\label{sec:examples:nmr}}
This test problem from magnetic resonance spectroscopy~\cite[Tab.~1]{VBH97}
considers sum of eleven complex exponentials.
\begin{equation}\label{eq:nmr}
	h_k = \sum_{k=1}^{11} a_k e^{135i \pi / 180} e^{(2i\pi f_k - d_k) j \delta}
\end{equation}
where $\delta = \frac13 \times 10^{-3}$ and parameters 
\begin{equation}
	\begin{matrix} 
	\ve a &= [& 75 & 150 & 75 & 150 & 150 & 150 & 150 & 150 & 1400 & 60 & 500 &]\phantom{.} \\
	\ve f &= [& -86&-70&-54&152&168&292&308&360&440&490&530 &]\phantom{.}\\
	\ve d &= [& 50&50&50&50&50&50&50&25&285.7&25&200 &].
	\end{matrix}
\end{equation}
To simulate detector noise, 
we add complex circular Gaussian random noise 
$\ve g\sim\mathcal{CN}(\ve 0, 15^2\; \ma I)$.
This example we use a total of $N=256$ measurements.

\begin{figure}
\centering 
\noindent
\begin{tabular}{@{}cc@{}}
\begin{tikzpicture}
	\begin{semilogyaxis}[
		ymin = 100,
		ymax = 2e4,
		ytickten = {2, 3, 4},
		height = 6.5cm,
		width = 8cm,
		xlabel = {$k$th singular value},
		ylabel = {$\sigma_k(\tma H)$},
		xtick = {1,2,...,20},
		xmin = 0.5, 
		xmax = 20.5,
		]
		\foreach \fn in {00, 01, 02, 03, 04, 05, 06, 07, 08, 09, 10, ..., 19}{ 
			\addplot[black, fill, fill opacity = 0.5] table [x=x,y expr=(10)^(\thisrow{y})] {fig_nmr_density_\fn.dat};
			\addplot[black, thick] table [x=x,y expr = (10)^(\thisrow{y})] {fig_nmr_extrema_\fn.dat};
		}
		\addplot[colorbrewerA1, thick] coordinates {(0.5,10^2.9958)  (20.5, 10^2.9958)}
			 [anchor = north west] node[pos=0] {upper bound (95\%)};
		\addplot[colorbrewerA2, thick] coordinates {(0.5,10^2.647)  (20.5, 10^2.647)}
			 [anchor = north west] node[pos=0] {empirical bound (95\%)};
	\end{semilogyaxis}
\end{tikzpicture}&
\begin{tikzpicture}
	\begin{groupplot}[
		group style = { group size = 1 by 3, vertical sep = 1.1cm},
		xmin = 8,
		xmax = 20,
		height = 2.5cm,
		width = 5.8cm,
		xtick = {1,2,...,20},
		axis y line = none,
		axis x line = bottom,
		clip = true,
	]
		\nextgroupplot[ymin = 0, ymax = 100, ylabel = frequency, ybar, clip = false]
		\addplot[const plot mark mid, fill = colorbrewerA1, fill opacity = 0.5, colorbrewerA1, 
			nodes near coords = {\pgfmathfloatifflags{\pgfplotspointmeta}{0}{}{\pgfmathprintnumber{\pgfplotspointmeta}\%}},
			every node near coord/.append style={anchor = west, rotate = 90, font=\tiny, yshift = -5pt, fill opacity = 1, xshift=0pt}] table
			 [x = rank, y expr = 100*\thisrow{bound}] {fig_nmr_count.dat};
		\addplot[black, const plot mark mid] coordinates {(10,0) (11,100) (12,0)};
		\node at (rel axis cs:1,1) [anchor=south east] {\cref{thm:order}}; 

		\nextgroupplot[ymin = 0, ymax = 100, ylabel = frequency, yshift = 0.0em, clip = false]
		\addplot[const plot mark mid, fill = colorbrewerA2, fill opacity = 0.5, colorbrewerA2, clip = true,
			nodes near coords = {\pgfmathfloatifflags{\pgfplotspointmeta}{0}{}{\pgfmathprintnumber{\pgfplotspointmeta}\%}},
			every node near coord/.append style={anchor = west, rotate = 45, font=\tiny, yshift = 5pt, fill opacity = 1, xshift = -2pt}] table
			[x = rank, y expr = 100*\thisrow{empirical}] {fig_nmr_count.dat};
		\addplot[black, const plot mark mid] coordinates {(10,0) (11,100) (12,0)};
		\node at (rel axis cs:1, 1) [anchor=south east] {Empirical Bound}; 

		\nextgroupplot[ymin = 0, ymax = 100, xlabel = estimated degree,ylabel = frequency, yshift = 0.0em, clip = false]
		\addplot[const plot mark mid, fill = colorbrewerA3, fill opacity = 0.5, colorbrewerA3, 
			nodes near coords = {\pgfmathfloatifflags{\pgfplotspointmeta}{0}{}{\pgfmathprintnumber{\pgfplotspointmeta}\%}},
			every node near coord/.append style={anchor = west, rotate = 45, font=\tiny, yshift = 5pt, fill opacity = 1, xshift = -2pt}] table
			[x=rank, y expr = 100*\thisrow{AIC}] {fig_nmr_count.dat};
		\addplot[black, const plot mark mid] coordinates {(10,0) (11,100) (12,0)};
		\node at (rel axis cs:1,1) [anchor=south east] {AIC};
	\end{groupplot}
\end{tikzpicture}
\end{tabular}

\caption{The application of our bounds and the AIC to estimate the number of complex exponentials
	embedded in complex Gaussian noise as described in \cref{sec:examples:nmr}.
	The left plot shows the distribution of the first twenty singular values of $\tma H_{129,128}$
	constructed from 1000 realizations, 
	where the frequency is denoted by the width of the shaded region
	and the range is denoted by the vertical black bar.
	The three right plots show the estimated model order using different techniques
	and the true model order of eleven is denoted by the hollow black rectangle.
}
\label{fig:nmr}
\end{figure}
 
\Cref{fig:nmr} illustrates how different approaches
estimate the number complex exponentials are present in $\tve h_N$;
recall from~\cref{sec:intro:related}, determining the number of exponentials
is equivalent to determining the McMillan degree.
As expected, the McMillan degree lower bound in~\cref{thm:order}
provides a lower bound on McMillan degree.
By using an empirical estimate of $\|\ma G_{m,n}\|_2$
as described in~\ref{sec:empirical}
we obtain a sharper and frequently correct
estimate of the McMillan degree.
This suggests that most of the loss of accuracy 
in our bound occurs mainly in the embedding step (\cref{lem:embedding}),
not in the use Weyl's theorem for the singular values \cref{thm:order}.
The AIC performs well in this case,
but as it requires identifying
requires more computation to obtain the fits for each candidate 
number of exponentials.

\subsection{Real Valued Data\label{sec:examples:real}}
As a second example, we consider the clamped beam model from the SLICOT benchmarks for model reduction~\cite{CD02}.
This considers a beam where the input is a force applied at the free boundary and the output is the displacement at this boundary.
Although originally a continuous time model, we can convert this to a discrete-time system in the form of~\cref{eq:ss}
using the matrix exponential
\begin{equation}\label{eq:beam}
\begin{split}
	\ve x_{j} &= e^{\ma A \delta} \ve x_{j-1} + \ve b u_j, \quad \text{where} \quad \ve b, \ve c,  \ve x_j \in \R^{348}, \ \ma A \in \R^{348 \times 348} \\
	y_j &= \ve c^\trans \ve x_{j}.
\end{split}
\end{equation}
Here  we take the time step $\delta=0.1$ and use $N=2^{13}=8192$ measurements
to which we add real Gaussian noise with $\ve g_N\sim \mathcal{N}(\ve 0, 10^{-2}\ma I)$.
Although this example has a McMillan degree of $348$, corresponding to the dimension of $\ma A$,
it is highly reducible and the singular values of $\ma H_{m,n}$ decay rapidly.
This simulates real systems which may have components that cannot be resolved due to noise.

\begin{figure}
\centering 
\noindent
\begin{tabular}{@{}cc@{}}
\begin{tikzpicture}
	\begin{semilogyaxis}[
		xmin = 0.5, 
		xmax = 20.5,
		ymin = 7,
		ymax = 2e4,
		height = 6.5cm,
		width = 8cm,
		xlabel = {$k$th singular value},
		ylabel = {$\sigma_k(\tma H)$},
		xtick = {1,2,...,20},
		]
		\foreach \fn in {00, 01, 02, 03, 04, 05, 06, 07, 08, 09, 10, ..., 19}{ 
			\addplot[black, fill, fill opacity = 0.5] table [x=x,y expr= (10)^(\thisrow{y})] 
				{fig_beam_density_\fn.dat};
			\addplot[black, thick] table [x=x,y expr = (10)^(\thisrow{y})] 
				{fig_beam_extrema_\fn.dat};
		}
		\addplot[colorbrewerA1, thick] coordinates {(0.5,10^1.63409711275)  (20.5, 10^1.63409711275)}
			 [anchor = north west] node[pos=0] {upper bound (95\%)};
		\addplot[colorbrewerA2, thick] coordinates {(0.5,10^1.27241427984)  (20.5, 10^1.27241427984)}
			 [anchor = north west] node[pos=0] {empirical bound (95\%)};
	\end{semilogyaxis}
\end{tikzpicture}&
\begin{tikzpicture}
	\begin{groupplot}[
		group style = { group size = 1 by 3, vertical sep = 1.1cm},
		xmin = 7,
		xmax = 20,
		height = 2.5cm,
		width = 5.8cm,
		xtick = {1,2,...,20},
		axis y line = none,
		axis x line = bottom,
	]
		\nextgroupplot[ymin = 0, ymax = 100, ylabel = frequency, ybar, clip = false]
		\addplot[const plot mark mid, fill = colorbrewerA1, fill opacity = 0.5, colorbrewerA1, 
			nodes near coords = {\pgfmathfloatifflags{\pgfplotspointmeta}{0}{}{\pgfmathprintnumber{\pgfplotspointmeta}\%}},
			every node near coord/.append style={anchor = south, rotate = 0, font=\tiny, xshift=0pt, yshift = -2pt, fill opacity = 1}] table
			 [x = rank, y expr = 100*\thisrow{bound}] {fig_beam_count.dat};
		\node at (rel axis cs:1,1) [anchor=south east] {\cref{thm:order}}; 

		\nextgroupplot[ymin = 0, ymax = 100, ylabel = frequency, yshift = 0.0em, clip = false]
		\addplot[const plot mark mid, fill = colorbrewerA2, fill opacity = 0.5, colorbrewerA2, clip = true,
			nodes near coords = {\pgfmathfloatifflags{\pgfplotspointmeta}{0}{}{\pgfmathprintnumber{\pgfplotspointmeta}\%}},
			every node near coord/.append style={anchor = south, rotate = 0, font=\tiny, yshift = -2pt, fill opacity = 1, xshift = 0pt}] table
			[x = rank, y expr = 100*\thisrow{empirical}] {fig_beam_count.dat};
		\node at (rel axis cs:1,1) [anchor=south east] {Empirical Bound}; 

		\nextgroupplot[ymin = 0, ymax = 100, xlabel = estimated degree,ylabel = frequency, yshift = 0.0em, clip = false]
		\addplot[const plot mark mid, fill = colorbrewerA3, fill opacity = 0.5, colorbrewerA3, 
			nodes near coords = {\pgfmathfloatifflags{\pgfplotspointmeta}{0}{}{\pgfmathprintnumber{\pgfplotspointmeta}\%}},
			every node near coord/.append style={anchor = west, rotate = 90, font=\tiny, yshift = 0pt, fill opacity = 1, xshift = 0pt}] table
			[x=rank, y expr = 100*\thisrow{AIC}] {fig_beam_count.dat};
		\node at (rel axis cs:1,1) [anchor=south east] {AIC};
	\end{groupplot}
\end{tikzpicture}
\end{tabular}

\caption{The application of our bounds and the AIC to estimate the McMillan degree
	of the beam model described in \cref{sec:examples:real}.
	The left plot shows the distribution of singular values of $\tma H_{4097,4096}$
	over one thousand realizations as in \cref{fig:nmr}.
	The three plots on the right show the estimated model order using different techniques.
}
\label{fig:beam}
\end{figure}

\Cref{fig:beam} illustrates different approaches for estimating the McMillan degree of this system.
Unlike the previous example, we have no hope of estimating the true McMillan degree of $\ma A$,
as even in the absence of noise only $105$ singular values of $\ma H_{4097,4096}$ exceed $10^{-10}$.
With the addition of noise 
we obtain a McMillan degree lower bound of $8$ using \cref{thm:order} and $12$ using the empirical estimate.
Both of these are smaller than than the McMillan degree estimate provided by the AIC.

 \section{Conclusion\label{sec:conclusion}}
Here we have established an upper bound on the norm 
of a random Hankel matrix with no unknown constants in~\cref{thm:main}
and used this result to construct a lower bound on the McMillan degree 
from noisy impulse response measurements in \cref{thm:order}.
As the examples in \cref{sec:examples} illustrate, 
this bound provides a useful lower bound on the McMillan degree
that can be applied to both modal analysis and system identification.
However in engineering practice, we expect the empirically determined 
bound on $\|\ma G_{m,n}\|_2$ to be more useful.
It provides a sharper bound and is easy to compute without 
knowledge of the underlying noise distribution
by using measurements of the system with no input.

 \section*{Acknowledgements}
I would like to thank Mark Embree for his support during my PhD
where this result originated,
Paul Martin for his feedback on an early draft of this manuscript,
and the anonymous reviewers for their help refining this manuscript.
 
\bibliographystyle{siamplain}
\bibliography{abbrevjournals,master}

\end{document}